\documentclass[fr,12pt]{article}
\usepackage[latin1]{inputenc}
\usepackage{amsmath,amsfonts,graphicx}

\def \R{\mathbb{R}} 
\def \N{\mathbb{N}}

\def \Z{\mathbb{Z}}
\def \N{\mathbb{N}}

\def \T{\mathbb{T}}

\def \AA{{\cal{A}}}

\def \GG{{\cal{G}}}
\def \GGb{\ov{{\cal{G}}}}
\def \GGe{{{\cal{G}}}^*}
\def \GGeb{\ov{\GG}^*}

\def \RR{{\cal{R}}}

\def \TT{{\cal{T}}}
\def \UU{{\cal{U}}}
\def \VV{{\cal{V}}}

\def \a{\alpha}
\def \b{\beta}
\def \g{\gamma}

\def \d{\delta}
\def \e{\varepsilon}
\def \g{\gamma}

\def \o{\omega}

\def \ov{\overline}

\def \Ab{\overline{A}}
\def \Bb{\overline{B}}
\def \dun{\partial_1}
\def \dde{\partial_2}
\def \duu{\partial_{11}}
\def \ddd{\partial_{22}}
\def \dud{\partial_{12}}

\def \eqv{\Longleftrightarrow}

\def \ft{\widetilde{\phi}}
\def \Fb{\ov{F}}
\def \Ker{{\rm Ker}}
\def \ob{\overline{\o}}
\def \Pb{\ov{\Phi}}
\def \pr{{\rm pr}}
\def \psib {\ov{\psi}}

\def \vv{\vert \vert}
\def \xd{\dot{x}}

\begin{document}

\bigskip
\bigskip
\centerline{\Large NON-RESONANT TORI IN SYMPLECTIC}
\centerline{\Large TWIST MAPS WITHOUT CONJUGATE POINTS}

\medskip
\centerline{Marc Arcostanzo}
\centerline{(Avignon University, LMA EA 2151, Avignon, France)} 

\bigskip
\bigskip

Let $d \geq 1$ be an integer, $\T^d$ the d-dimensional torus, and 
$$F : T^*\T^d \longrightarrow T^*\T^d$$
\noindent 
a $C^\infty$ twist map. Twist maps are ewamples of symplectic diffeomorphisms ; they will be defined 
 more precisely in part 1. We assume that $F$ is without conjugate points. This means that 
$$ \forall n \in \Z \setminus \{ 0 \}, \forall (\xd,p) \in T^*\T^d, 
DF^n(\VV(\xd,p)) \cap \VV(F^n(\xd,p)) = \{ 0 \} ,$$   
\noindent
where $\VV(\xd,p)$ denotes the vertical space at the point $(\xd,p)$. It is proved in \cite{Arc} that 
this implies that $F$ is $C^0$-integrable, i.e. there exists a continuous foliation of $T^*\T^d$, each 
leaf being a Lipschitz Lagrangian graph that is $F$-invariant. 

\smallskip
It could be that all the leaves of this foliation are in fact smooth, but this is still an open question. 
However, some of them are indeed smooth : it is shown in \cite{Arc} (proposition 3.1) that if $\Fb$ is 
a lift of $F$ to $T^*\R^d$ then for every $n \in \N$ and every $r \in \Z^d$, the set 
$$\GGeb_{N,r} = \{ (x,p) \in T^*\R^d \ {\rm s.t.} \ \Fb^N(x,p) = (x+r,p) \}$$
\noindent
is a $C^\infty$ Lagrangian $\Fb$-invariant graph. Its projection $\GGe_{N,r}$ on $T^*\T^d$ is one 
of the leaves of the foliation. It is by definition a union of periodic orbits sharing the same 
period $N$.

\smallskip
Here we study the dynamics of $F$ in a neighborhood of $\GGe_{N,r}$. We use a KAM theorem to show 
the existence of a rich family of $F$-invariant Lagrangian graphs accumulating on $\GGb^*_{N,r}$, 
on which $F$ is conjugated to a translation of non-resonant vector. In fact, we have 
the following result : 

\bigskip
{\bf Theorem } : {\it Let $\ob \in \R^d$ be strongly Diophantine vector, i.e. there are real numbers 
$\g > 0$ and $\tau > 0$ such that 
$$ \forall k \in \Z^d \setminus \{ 0 \}, \ \forall l \in \Z, \ 
\vert k \cdot \ob + l \vert \geq {\g \over \vert k \vert ^\tau}.$$

For every large integer $m$, there is a $C^\infty$ 
Lagrangian embedding $i_m : \T^d \longrightarrow T^*\T^d$ such that 

i) $i_m(\xd) = (\psi_m(\xd), f_m(\xd))$, where $\psi_m$ is a $C^\infty$ diffeomorphism of $\T^d$, isotopic 
to the identity, and 
$$ \TT_m = i_m(\T^d) = \{ (\xd , (f_m \circ \psi_m^{-1})(\xd)) ; \xd \in \T^d \} $$
\noindent
is a Lagrangian graph ; the sequence $(\TT_m)$ converges to $\GGe_{N,r}$ in $C^\infty$ topology.

ii) The sequence $(\psi_m)$ converges in $C^\infty$ topology to a diffeomorphism $\psi_\infty$ (independant 
of $\ob$), isotopic to the identity.

iii) The tori $\TT_m$ are $F$-invariant and the restriction of $F$ to $\TT_m$ is conjugated to a non-resonant 
translation. More precisely,
$$ \forall n \in \Z, \ \forall \xd \in \T^d, \  F^n(i_m(\xd)) = i_m(\xd + {n \over N}r + {n \over mN}\ob).$$} 

Note that this gives us some insight into the dynamics of $F$ restricted to $\GGe_{N,r}$.

\bigskip
{\bf Corollary } : {\it The diffeomorphism $\psi_\infty$ conjugates the action of $F$ on $\GGe_{N,r}$ to a 
translation of vector ${r\over N}$ on $\T^d$}.

\bigskip
In the case of a continuous flow associated to a Tonelli Hamiltonian, a similar result is established in 
\cite{AABZ}. Our strategy is to mimic the proof given in this article, but many problems arise when we switch 
from the continuous to the discrete case. For example, we can no more make use of a quantity that is constant 
along the orbits (as the Hamiltonian in the continous case), or derivate along the flow. As a result, some 
key parts of the proof need totally different arguments. 

\smallskip
The paper is organized as follows. In section 1, we briefly recall some basic facts on twist mapsand some results 
of \cite{Arc}. In section 2, we explain how to find a normal form for $F$ in the neighborhood of $\GGb^*_{N,r}$. 
This requires two lemmas which are proved in section 3. We then apply a KAM theorem in section 4 and  
explain the end of the proof of the theorem. 
  
\bigskip
\bigskip
{\Large 1. Twist maps without conjugate points}

\bigskip
\medskip
Here we give a brief introduction to the theory of twist maps. We refer the reader to \cite{G} for a complete study.
Let $d \geq 1$ be an integer. Denote by $\T^d = \R^d / \Z^d$ 
the d-dimensional torus. Let $\T^*\R^d = \R^d \times (\R^d)^*$ be the cotangent space of $\R^d$. Consider a generating 
function, that is a map $S : \R^d \times \R^d \longrightarrow \R$ of class $C^\infty$ which satisfies the following 
two conditions :  

\noindent
$(C1) \ \forall r \in \Z^d, \forall (x,y) \in \R^d \times \R^d, S(x+r,y+r) = S(x,y)$ ;

\noindent
$(C2)$ (`uniform twist condition', see \cite{BMK}) There is a real number $A > 0$ for which  

$$ \forall (x,y) \in \R^d \times \R^d, \forall \xi \in \R^d, \ 
\sum _{i,j} {\partial^2S(x,y) \over \partial x_i \partial y_j}(x,y) \xi_i \xi_j \leq - A \vv \xi \vv^2. $$

\smallskip
A sequence $(x_n)_{n \in \Z}$ with values in $\R^d$ is said to be extremal if it satisfies 
$$\forall n \in \Z, \ \dde S(x_{n-1} , x_n) + \dun S(x_n , x_{n+1}) = 0.$$
\noindent
Extremal sequences are the critical points of the (formal) action functionnal which assigns to 
each sequence $(x_n)_{n \in \Z}$ the sum of the serie $\sum_{n \in \Z} S(x_n,x_{n+1})$. The 
generating function also gives rise to a symplectic diffeomorphism $F$ of $T^*\T^d$. Let 
$\Fb : T^*\R^d \longrightarrow T^*\R^d$ be the diffeomorphism implicitely defined by  
$$ \Fb(x,p) = (x',p') \eqv p = -\dun S (x , x') \ {\rm and} \ p' = \dde S (x , x').$$
\noindent 
The diffeomorphism $\Fb$ is exact symplectic, which means that $\Fb^* \a - \a = dS$, where 
$\a = \sum_{i=1}^d x_i dq_i$ is the Liouville 1-form on $T^*\R^d$. Note that condition $(C1)$ 
implies that $\Fb$ is the lift to $\TT^*$ of a symplectic diffeomorphism $F$ of $T^*\T^d$.

\smallskip
Let $\pr_1 : (\xd,p) \in T^*\T^d \longmapsto \xd \in \T^d$ be the canonical projection. The 
vertical space at $(\xd,p) \in T^*\T^d$ is $\VV(\xd,p) = \Ker D\pr_1 (\xd,p)$. We say that the twist 
map $F$ is without conjugate points if  
$$ \forall n \in \Z \setminus \{ 0 \}, \ \forall (\xd,p) \in T^*\T^d,  \ 
DF^n(\VV(\xd,p)) \cap \VV(F^n(\xd,p)) \ = \ \{ 0 \}.$$
\noindent
This hypothesis has strong consequences on the behaviour of extremal sequences. It is shown in \cite{Arc} (corollary 1.5) 
that if $F$ is without conjugate points, then for every $(x,y) \in \R^d \times \R^d$ and every integer 
$N \geq 1$, there is a unique extremal sequence $(x_n)_{n \in \Z}$ with $x_0 = x$ et $x_N = y$. Moreover, 
this extremal sequence minimizes the action in the following sense. If $k$ and $l$ are two integers 
with $l - k \geq 2$, then, letting $y_k = x_k$ and $y_l = x_l$, one has 
$$\forall (y_{k+1}, \dots , y_{l-1}) \in (\R^d)^{l-k-1}, \ 
\sum_{i=k}^{l-1} S(x_i, x_{i+1}) \leq \sum_{i=k}^{l-1} S(y_i, y_{i+1}),$$
\noindent
and equality holds if and only if $y_i = x_i$ for every $i \in \{ k+1, \dots, l-1 \}$. 

\smallskip
Let $r \in \Z^d$, $N \geq 1$ an integer, and $x \in \R^d$. Consider the extremal sequence $(x_n)_{n \in \Z}$ 
with $x_0 = x$ and $x_N = x_0 + r$. It is a non-trivial fact (see \cite{Arc}, proposition 2.1) that we have
$$ \forall n \in \Z, \ x_{n+N} = x_n +r.$$ 
\noindent
We can use this to construct periodic orbits of $F$. In fact, define a sequence $(p_n)$ with values 
in $(\R^d)^*$ as follows : 
$$\forall n \in \Z, \ p_n = -\dun S (x_n , x_{n+1}).$$
\noindent
It is periodic (this is a consequence of (C1)). Morever, $(x_n,p_n)_{n \in \Z}$ is an orbit of $\Fb$ whose 
projection to $T^*\T^d$ is periodic with period $N$. Letting $x$ vary in $\R^d$, we obtain a subset 
$\GGe_{N,r}$ of $T^*\T^d$ that is a graph over the whole of $\T^d$. It is a union of periodic orbits of $F$ 
and is therefore $F$-invariant. A result of \cite{Arc} (proposition 3.1) is that this graph is of class 
$C^\infty$ and Lagrangian.   

The aim of this paper is to study the dynamics of $F$ in a neighborhood of $\GGe_{N,r}$. We will make use of 
the following tool. Assume that $\Phi$ is a symplectic diffeomorphism of $T^*\T^d$ that leaves invariant the 
null section $0_{\T^d} = \{ (\xd,0), \xd \in \T^d \}$. For $\e > 0$, consider the map 
$$\RR_\e : (x,p) \in T^*\T^d \longmapsto (x,\e p) \in T^*\T^d.$$
\noindent
Since $\RR_\e^*\a = \e \a$, $\Phi_\e = \RR_\e^{-1} \circ \Phi \circ \RR_\e$ is, just as $\Phi$,  a symplectic 
diffeomorphism of $T^*\T^d$ with $\Phi_\e(0_{\T^d}) = 0_{\T^d}$. The study of $\Phi_\e$ when $\e$ goes to 
$0$ gives us an insight into how $\Phi$ behaves near $0_{\T^d}$. Note that when $\Phi =F$ is the twist map 
associated to the generating function $S$, then $\Phi_\e$ is the twist map associated to the generating 
function ${S \over \e}$.

We shall use the following notations. $x$ always denotes a point in $\R^d$, while $\xd$ refers to an element 
of $\T^d$. $< \cdot , \cdot > : \R^d \times (\R^d)^* \longrightarrow \R$ is the duality bracket. 
If $M$ is a matrix or a linear operator, we note $M^T$ its transpose and (if $M$ is invertible) $M^{-T}$ the 
inverse of $M^T$. 

\bigskip
\bigskip
{\Large 2. A normal form for $F^N$}

\bigskip
\medskip
We fix once and for all an integer $N \geq 1$ and $r \in \Z^d$. The Lagrangian graph $\GGe_{N,r}$ may be 
written as
$$ \GGe_{N,r} = \{ (\xd, p_\infty + du(\xd)), \xd \in \T^d \},$$
\noindent
with $p_\infty \in (\R^d)^*$ and $u : \T^d \longrightarrow \R$ a $C^\infty$ map. In this section, we explain 
how to obtain a normal form for $F^N$ in the neighborhood of $\GGe_{N,r}$. The precise statement is as 
follows.

\medskip
{\bf Proposition 1} : {\it There exists a symplectic $C^\infty$ diffeomorphism $G$ of $T^*\T^d$ of the form
$$ G(\xd,p) = (\psi(\xd), p_\infty + Du(\psi(\xd)) + D\psi(\xd)^{-T}p), $$
\noindent
where $\psi$ is a diffeomorphism of $\T^d$ isotopic to the identity, such that $G(0_{\T^d}) = \GGe_{N,r}$, 
and 
$$G^{-1} \circ F^N \circ G(\xd,p) = (\xd + \Bb p + O(p^2), p + O(p^3)),$$
\noindent
where $\Bb \in L((\R^d)^*,\R^d)$ is symmetric positive definite.}

\medskip
Proof : Consider the symplectic change 
of variables 
$$ G_0(\xd,p) = (\xd, p + p_\infty + du(\xd)), $$
and the generating function $S_0 (x,y) = S(x,y) - u(x) + u(y) + < p_\infty , y - x >.$ It satifies the conditions 
(C1) and (C2). The associated exact symplectic diffeomorphism of $T^*\T^d$ is $F_0 = G_0^{-1} \circ F \circ G_0$. 
Since $G_0$ preserves the fibers and its restriction to each fiber is a translation, $F_0$ is without conjugate 
points. By definition of $\GGb^*_{N,r}$ we have 
$$\forall \xd \in \T^d, \ F_0^N(\xd,0) = (\xd,0).$$
This implies that the differential of $F_0^N$ at $(\xd,0)$ takes the form
$$ D F_0^N (\xd,0) [\d x,\d p] = (\d x + B(\xd)\d p , D(\xd)\d p), $$
\noindent
with $B(\xd) \in L((\R^d)^*,\R^d)$ and $D(\xd) \in L((\R^d)^*,(\R^d)^*)$. Moreover $D F_0^N (\xd,0)$ is a 
symplectic linear map, so $D(\xd) = id_{(\R^d)^*}$ and $B(\xd)$ is symmetric. In fact, we can say more about 
$B(\xd)$ : 

\medskip
{\bf Lemma 1} : {\it $B(\xd)$ is symmetric positive definite.}

\medskip 
This lemma will be proved in the next section. As a consequence, we can define a Riemannian metric $g$ on $\T^d$ : 
$$\forall \xd \in \T^d, \ \forall v,v' \in \R^d, \ g((\xd,v),(\xd,v')) = <B(\xd)^{-1}v,v'>.$$
\noindent
The next step is to prove that $g$ enjoys a rather strong property : 

\medskip
{\bf Lemma 2} : {\it The metric $g$ is without conjugate points.}

\medskip
Once again, we postpone the proof to the next section. D. Burago and S. Ivanov proved (see \cite{BI}) that any 
Riemannian metric on the torus that is free of conjugate points must be flat. So $g$ is flat : there 
exists a $C^\infty$ diffeomorphism $\psi$ of $\T^d$ isotopic to the identity and a symmetric positive definite 
$\Ab \in L_(\R^d,(\R^d)^*)$ such that   
$$\forall \xd \in \T^d, \ \forall v,v' \in \R^d, 
<B(\psi(\xd))^{-1}D\psi(\xd) \cdot v, D\psi(\xd) \cdot v'> = <\Ab v,v'>.$$
\noindent
Let $\Bb = \Ab^{-1}$. Then 
$$ \forall \xd \in \T^d, \ \Bb = D\psi^{-1}(\xd) B(\psi(\xd))D\psi(\xd)^{-T}.$$
\noindent
Consider the symplectic diffeomorphism of $\T^d$ : $G_1(\xd,p) = (\psi(\xd), D\psi(\xd)^{-T})$ and the 
composition $G = G_0 \circ G_1$. We have 
$$G(\xd,p) = (\psi(\xd),  p_\infty + Du(\psi(\xd)) + D\psi(\xd)^{-T}p).$$
Let $F_1 = G_1^{-1} \circ F_0 \circ G_1  = G^{-1} \circ F \circ G$. It satisfies 
$$\forall \xd \in \T^d, \ F_1^N (\xd,0) = (\xd,0) \ {\rm and} \ D_pF_1^N(\xd,0) = \Bb.$$
\noindent 
Then the next lemma applied to $F_1^N$ implies that $F_1^N(\xd,p) = (\xd + \Ab p + O(p^2) , p + O(p^3))$ 
as desired. 

\medskip
{\bf Lemma 3} : {\it Let $F: T^*\T^d \longrightarrow T^*\T^d$ be an exact symplectic diffeomorphism. 
Assume that $F$ fixes every point of $0_{\T^d}$. Let $\Bb(\xd) = D_pF(\xd,0)$ for every $\xd \in \T^d$. 
Then 
$$F(\xd,p) = (\xd + \Bb(\xd) p + O(p^2) , p - {1 \over 2} D_p<B(\xd)p,p> + O(p^3)).$$}

\medskip
The proof of this lemma is given in \cite{AABZ} (page 182). Note that it uses in a crucial way 
the maps $\Phi_\e$ introduced at the end of the first section.

\vfill \eject
{\Large 3. Proof of lemma 1 and lemma 2}

\bigskip
\medskip

We begin with the proof of lemma 1. Let $\xd_0 \in \T^d$ and $(\xd_n)$ the sequence of points in $\T^d$ 
defined as $F_0^n(\xd_0,0)=(\xd_n,0)$. As $F_0$ preserve $0_{\T^d}$, the matrix of $DF_0(\xd_i,0)$ in the 
canonical basis is a symplectic matrix of the form
$$ m_i = \begin{pmatrix} a_i & b_i \\ O & d_i \end{pmatrix},$$
\noindent
hence the matrix $s_i := ^tb_i d_i$ symmetric. It is straightforward to check that if two matrices of the 
type $m_i$ are such that the $s_i$ are positive definite, then the same property holds for their product. 
By an immediate induction, the same is true for a product of any number of such matrices.

\smallskip
For all $\xd \in \T^d$, $F_0^N(\xd,0)=(\xd,0)$, so that the matrix of $DF_0^N(\xd_0,0)$ may be written as
$$ M_N = \begin{pmatrix} I_d & B_N \\ O_d & D_N \end{pmatrix}.$$
\noindent
$M$ is symplectic, hence $D_N = I_d$ and $B_N$ is symmetric. As a consequence of the chain rule, 
$M = m_{N-1}m_{N-2} \dots m_1m_0$. So it only remains to prove that each $s_i$ is positive definite 
to get the conclusion that $B_N$ is also positive definite. As $\xd_0$ is an arbitrary point in $\T^d$, 
we only need to check that $s_0$ is positive definite.  

\smallskip
It is possible to express the differential of a twist map in terms of its generating function : the 
result is that
$$ m_0 = \begin{pmatrix} -\dud S(x_0,x_1)^{-1} \duu S(x_0,x_1) & -\dud S(x_0,x_1)^{-1} \\ 
O_d & -\ddd S(x_0,x_1) \dud S(x_0,x_1)^{-1} \end{pmatrix}.$$
\noindent
Therefore  $s_0 = \dud S(x_0,x_1)^{-T} \ddd S(x_0,x_1) \dud S(x_0,x_1)^{-1}$, and we have to show that 
$\ddd S(x_0,x_1)$ is positive definite to finish the proof.

\smallskip
Let $\Fb_0$ be a lift of $F_0$ to $T^*\R^d$ and $x_0 \in \R^d$. For any integer $n$, we have $\Fb_0^n(x_0,0) 
= (x_n,0)$ for some $x_n \in \R^d$. For every $p \in (\R^d)^*$, let $x_n(p) = \pr_1(\Fb_0^n(x_0,p))$. Clearly 
$x_0(p) = x_0$ for all $p$, and $x_n = x_n(0)$ for all $n$. Consider the action fonctional 
$$\AA_n(p) = S(x_0, x_1(p)) + \sum_{i=1}^{n-1}S(x_i(p),x_{i+1}(p)) \ + \ S(x_n(p),x_{n+1}).$$
Since $\Fb$ is without conjugate points, we know (see part 1) that $\AA_n$ admits a global strict 
minimum at $p=0$. We compute 
$$D\AA_n(p) = \dde S(x_0,x_1(p))Dx_1(p) + \sum_{i=1}^{n-1} D_i + \dun S(x_n(p),x_{n+1})Dx_{n+1}(p),$$
$${\rm where} \ \  D_i = \dun S(x_i(p),x_{i+1}(p)) Dx_i(p) + \dde S(x_i(p),x_{i+1}(p)) Dx_{i+1}(p).$$
As $(x_0,x_1(p),x_2(p), \dots,x_{n+1}(p))$ is an extremal sequence, we have 
$$\dde S(x_{i-1}(p),x_{i}(p)) + \dun S(x_i(p),x_{i+1}(p)) = 0$$ 
\noindent
for all $i$, so almost all terms in $D\AA_n(p)$ cancel out and we are left with 
$$D\AA_n(p) = [\dde S(x_{n-1}(p),x_n(p)) + \dun S(x_n(p),x_{n+1}(p))] Dx_n(p).$$
This may be rewritten as 
$$D\AA_n(p) = [-\dun S(x_n(p),x_{n+1}(p)) + \dun S(x_n(p),x_{n+1}(p))] Dx_n(p),$$
\noindent
whence the following expression for the second differential of $\AA_n$ : 
$$D^2\AA_n(0) : (v,v') \longmapsto - \dud S(x_n,x_{n+1}) (Dx_n(0) \cdot v, Dx_{n+1}(0) \cdot v') .$$
This implies that $S_n = - B_n^T \dud S(x_n,x_{n+1}) B_{n+1}$ is the (symmetric positive definite) 
matrix of $D^2\AA_n(0)$.

\smallskip
From now on, we assume that $n = kN$ for some integer $k$. Hence we have $x_n = x_0 + kr$, $x_{n+1} = 
x_1 + kr$ and condition (C1) implies $\dud S(x_n,x_{n+1}) = \dud S(x_0,x_1)$. Using the same argument 
and the chain rule, the matrix of $DF_0^n(x_0,0) = DF_0^{kN}(x_0,0)$ is 
$$ M_{kN} = \begin{pmatrix} I_d & B_{kN} \\ O_d & I_d \end{pmatrix}
 = \begin{pmatrix} I_d & B_{N} \\ O_d & I_d \end{pmatrix} ^k 
 = \begin{pmatrix} I_d & kB_{N} \\ O_d & I_d \end{pmatrix},$$
\noindent
so that $B_n = B_{kN} = kB_N$. 

\smallskip
To compute $B_{kN+1}$, we use once again the matrix $m_0$ introduced above. The relation $D\Fb_0^{kN+1}(x_0,0) =  
D\Fb_0(x_{kN},0) \circ D\Fb_0^{kN}(x_0,0)$ implies that $M_{kN+1} = m_0 M_{kN}$, so we get
$$B_{kN+1} = -k \dud S(x_0,x_1)^{-1} \duu S(x_0,x_1) B_{N} - \dud S(x_0,x_1)^{-1},$$
\noindent 
and finally
$$S_{kN} = k^2 B_N \duu S(x_0,x_1) B_N + k B_N.$$

Dividing by $k^2$ and letting $k$ go to infinity implies that $B_N \duu S(x_0,x_1) B_N$ is positive 
semi-definite. But $F_0$ has no conjugate points, so $B_N$ is invertible and $\duu S(x_0,x_1)$ is 
positive semi-definite. The matrix $m_0$ being symplectic, we have
$$\duu S(x_0,x_1) ^t \dud S(x_0,x_1)^{-1} \ddd S(x_0,x_1) \dud S(x_0,x_1)^{-1} = I_d.$$
\noindent
This implies that $\duu S(x_0,x_1)$ is invertible, so it has to be positive definite, as well as 
$\ddd S(x_0,x_1)$.

\medskip
The proof of lemma 2 is close to the one given \cite{AABZ} (page 182 and 183), so we will only sketch 
it, trying to put into perspective the main ideas, and referring to \cite{AABZ} for technical details. 
We lift $g$ to a $\Z^d$-periodic Riemannian metric on $\R^d$ and consider the corresponding Hamiltonian 
function  
$$H : (x,p) \in T^*\R^d \longmapsto {1 \over 2} <B(x)p,p> \in \R,$$ 
\noindent
with Hamiltonian vector field $X_H$ and Hamiltonian flow $(\ft_t^H)_{t \in \R}$. To prove the absence of conjugate 
points, we argue by contradiction. As explained in \cite{AABZ}, if $g$ had conjugate points then we could find two 
points $x$ and $y$ in $\R^d$ connected in time $S > 0$ by two distinct non-degenerate geodesics. So there would be 
$p_1$ and $p_2$ in $(\R^d)^*$ such that 
\begin{equation}   \pr_1 \circ \ft_S^H (x,p_1) = \pr_1 \circ \ft_S^H (x,p_2) = y \end{equation}
\noindent
with $D_p\pr_1 \circ \ft_S^H (x,p_1)$ and $D_p\pr_1 \circ \ft_S^H (x,p_2)$ invertible. 

We now consider $\Fb_0 : T^*\R^d \longrightarrow T^*\R^d$ a lift of $F_0$. If we could replace $\ft_S^H$ by $\Fb_0^N$ 
in $(1)$, we would get a contradiction since $F_0$ is without conjugate points. So we try to find some 
link between $\ft_S^H$ and $\Fb_0^N$. Once again, consider the maps 
$$\Phi_\e = \RR_\e^{-1} \circ F_0^N \circ \RR_\e : T^*\T^d \longrightarrow T^*\T^d$$
\noindent
and a lift $\Pb_\e : T^*\R^d \longrightarrow T^*\R^d$ (chosen in such 
a way that $\Pb_\e$ is close to the identity when $\e \rightarrow 0$). Then lemma 3 implies that  
$$\Pb_\e(x,p) = (x + \e B(x)p + O(\e^2), p - {\e \over 2} D_x <B(x)p,p> + O(\e^2))$$
\noindent 
and hence 
$$\Pb_\e(x,p) = (x,p) + \e X_H(x,p) + O(\e^2).$$

This is reminiscent of the Euler method used for numerical integration of ordinary differential equations. 
So we may hope that if $\e$ is small enough and $n$ not too large, $\Pb_\e^n(x,p)$ won't be very far from 
$\ft_{n\e}^H(x,p)$.  This turns out to be true, we refer to Lemma 2.1 of \cite{AABZ} for a rigorous statement. 
Applied to our case, it implies that when $m$ goes to infinity, the sequence $(\Pb_{S/m}^m)$ converges to 
$\ft_S^H$ on coimpact sets in topology $C^1$. As $D_p\pr_1 \circ \ft_S^H (x,p_1)$ and 
$D_p\pr_1 \circ \ft_S^H (x,p_2)$ are invertible, we may use the implicit function theorem to obtain that 
when $m$ is large enough, we can find $p'_1$ close to $p_1$ and $p'_2$ close to $p_2$ such that 
$$\pr_1 \circ \Pb_{S/m}^m (x,p'_1) = \pr_1 \circ \Pb_{S/m}^m (x,p'_2) = y.$$
\noindent
It is easy to check that we then have 
$$ \pr_1 \circ \Fb_0^{mN} (x,p'_1) = y + mr = \pr_1 \circ \Fb_0^{mN} (x,p'_2),$$
\noindent
and hence a contradition.

\vfill \eject
{\Large 4. Construction of the non-resonant tori}
\bigskip
\medskip

The existence of these tori is given by the following proposition. Its proof is not complicated and is very 
similar to the proof of Proposition 7 in \cite{AABZ}, so we will not repeat it here. The main idea is to apply 
a KAM theorem (theorem 1.2.3 in \cite{Bost}) to the family of symplectic symplectic maps $\UU_m = \Phi_{1/m}^m = 
\RR_{1/m}^{-1} \circ F^m \circ \RR_{1/m}$.  

\medskip
{\bf Proposition 2} : {\it Let $F : T^*\T^d \longrightarrow T^*\T^d$ be a $C^\infty$ symplectic diffeomorphism. 
Assume that on a neighborhood of $0_{\T^d}$, 
$$ F(\xd,p) = (\xd + \Bb p + O(p^2), p + O(p^3)),$$
\noindent
where $\Bb \in L((\R^d)^*,\R^d)$ is symmetric non-degenerate. Let $\ob \in \R^d$ be strongly Diophantine. 
Then for any large $m$ there is a $C^\infty$ Lagrangian embedding $j_m : \T^d \longrightarrow T^*\T^d$ such that 
$$ \forall \xd \in \T^d, F^m(j_m(\xd)) = j_m(\xd + \ob). $$
\noindent
Moreover, $j_m$ is of the following form : 
$$j_m(\xd) = (\xd + u_m(\xd), \Bb^{-1} \left ( {\ob \over m} \right ) + v_m(\xd)),$$
\noindent
with $u_m : \T^d \longrightarrow \R^d$ and $v_m : \T^d \longrightarrow (\R^d)^*$ of class $C^\infty$, and, 
for any $k$, 
$$\vv u_m \vv_{C^k(\T^d)} = o(1) \ {\rm and} \ \vv v_m \vv_{C^k(\T^d)} = o \left ( {1 \over m} \right ) 
\ {\rm as} \ m \rightarrow \infty .$$}

\smallskip
Let $F_1 = G^{-1} \circ F \circ G$, where $G$ is the symplectic diffeomorphism given by proposition 1. 
According to this proposition, we have 
$$ F_1^N(\xd,p) = G^{-1} \circ F^N \circ G(\xd,p) = (\xd + \Bb p + O(p^2), p + O(p^3)),$$
\noindent
so that we may apply proposition 2 to $F_1^N$. Consider the set $j_m(\T^d)$. It is clearly invariant 
by $F_1^{Nm}$, but the following stronger result holds.

\medskip
{\bf Lemma 4} : {\it $j_m(\T^d)$ is $F_1$-invariant.}

\medskip
{\it Proof} : $F$ is without conjugate points, so $T^*\T^d$ is the disjoint union of $F$-invariant graphs 
$(g_i)_{i \in I}$, and hence a disjoint union of the $F_1$-invariant graphs $(G^{-1}(g_i))_{ i\in I}$. Pick 
$j_m(\xd_0,p_0) \in j_m(\T^d)$, it belongs to some $G^{-1}(g_{i_0})$. As $G^{-1}(g_{i_0})$ is $F_1$-invariant, we only 
need to show that $j_m(\T^d) = G^{-1}(g_{i_0})$. Consider
$$ E = \{ \xd \in \T^d \ {\rm s.t.} \ j_m(\xd) \in G^{-1}(g_{i_0}) \}.$$
\noindent
Note that if $\xd \in E$, then so does $\xd + \ob$, since $j_m(\xd + \ob) = F_1^{Nm}(j_m(\xd))$ and $G^{-1}(g_{i_0})$ 
is $F_1$-invariant. The point $\xd_0$ belongs to $E$, hence $E$ contains $\xd_0 + k \ob$ for all integer $k$. These 
points are dense in $\T^d$ and $E$ is clearly closed, so $E = \T^d$ and $j_m(\T^d) \subset G^{-1}(g_{i_0})$. Then the 
compact manifold $j_m(\T^d)$ is included in the connected manifold $G^{-1}(g_{i_0})$, and they have the same dimension, 
so they coincide. \hfill{$\square$}

\smallskip
For any integer $n$ and any (large) integer $m$, we may then consider the map 
$$\a_{m,n} = j_m^{-1} \circ F_1^n \circ j_m : \T^d \longrightarrow \T^d.$$
\noindent
According to proposition 2, $\a_{m,Nm}$ is the translation $\tau_{\ob}$ of vector $\ob$. Since $F_1^n$ commutes 
with $F_1^{Nm}$, $\a_{n,m}$ commutes with $\a_{m,Nm} = \tau_{\ob}$. Using the same topological arguments as in the proof of 
lemma 4, we may conclude that every $\a_{m,n}$ is a translation. Moreover we clearly have $\a_{n+1,m} = \a_{n,m} 
\circ \a_{1,m}$ for every integer $n$, so if $\a_{m,1}$ is the translation of vector $\b_m$, then $\a_{m,n}$ is 
the translation $n\b_m$. As $\a_{m,Nm} = \tau_{\ob}$, we get
\begin{equation} \exists k_m \in \Z^d \ {s. t.} \ Nm\b_m = \ob + k_m.  \end{equation}

We are going to show that if $m$ is large enough, then $k_m = mr$. When $m$ goes to infinity, 
$j_m$ converges uniformly to $j_\infty : \xd \in \T^d \longmapsto (\xd,0) \in T^*\T^d$. 
Let $f_1 : \T^d \longrightarrow \T^d$ the map such that $F_1(\xd,0) = (f_1(\xd),0)$ for all $\xd \in \T^d$. Then 
$$\pr_1(j_m(\xd + \b_m)) = \pr_1 \circ j_m \circ \a_{m,1} (\xd) = \pr_1 \circ F_1 \circ j_m (\xd).$$
The right-hand side converges to $f_1(\xd)$, while the left-hand is equal to 
$$\xd + \b_m + u_m(\xd + \b_m) = \xd + \b_m + o(1).$$ 
\noindent
This implies that the sequence $(\b_m)$ converges to some vector $\b_\infty$ and that $f_1$ is the translation of 
vector $\b_\infty$.

\smallskip
Recall that the map $\psi : \T^d \longrightarrow \T^d$ given by proposition 1 is isotopic to the identity. So if 
$\psib$ denotes a lift of $\psi$ to $\R^d$, $\psib$ commutes with the translation of vector $r$. Let $\ov{G}$ be 
a lift of $G$ to $T^*\R^d$ with $\pr_1 \circ \ov{G} = \psib$. Then $\Fb_1 = \ov{G}^{-1} \circ \Fb \circ \ov{G}$ 
is a lift of $F_1$ to $T^*\R^d$. We now compute $\Fb_1^N(x,0) = \ov{G}^{-1} \circ \Fb^N \circ \ov{G} (x,0)$ for 
any $x \in \R^d$. To begin with, $\ov{G} (x,0) = (\psib(x),p)$ for some $p \in (\R^d)^*$. As $\ov{G}(0_{\R^d}) = 
\GGeb_{N,r}$, $\Fb^N (\psib(x),p) = (\psib(x)+r,p)$. Since$\psib$ commutes with the translation of vector $r$, we 
finally have $\Fb_1^N(x,0) = \ov{G}^{-1} (\psib(x)+r,p) = (x+r,0)$.  Since $F_1(\xd,0) = (\xd + \b_\infty,0)$, 
this implies that $\b_\infty = {r \over N}$, so that
\begin{equation}
N \b_m = r + o(1) \end{equation}

We know that $\a_{m,N} = j_m^{-1} \circ F_1^N \circ j_m$ is the translation of vector $N\b_m$. This implies
$$ \forall \xd \in \T^d, \ F_1^N( j_m(\xd)) = j_m(\xd + N\b_m). $$
\noindent
According to proposition 2, we have 
$$F_1^N( j_m(\xd)) = F_1^N( \xd + u_m(\xd) , \Bb^{-1} \left ( {\ob \over m} \right )  + v_m(\xd)).$$
\noindent
Using the estimates on $u_m$ and $v_m$ given by proposition 2 and proposition 1, we get
\begin{equation}
F_1^N( j_m(\xd)) = (\xd + u_m(\xd) + {\ob \over m} + o \left ( {1 \over m} \right ) , 
\Bb^{-1} \left ( {\ob \over m} \right )  + o \left ( {1 \over m} \right ) ) \end{equation}
\noindent
On the other hand, 
\begin{equation} 
j_m(\xd + N\b_m) = (\xd + N\b_m + u_m(\xd + N\b_m) , \Bb^{-1} \left ( {\ob \over m} \right ) 
+ o \left ( {1 \over m} \right ) )  \end{equation}
\noindent
Comparing $(4)$ and $(5)$ leads to : there is a vector $l_m \in \Z^d$ such that 
$$\forall \xd \in \T^d, \ N\b_m + u_m(\xd + N\b_m) = u_m(\xd) + {\ob \over m} + l_m + 
o \left ( {1 \over m} \right ) $$
\noindent
Taking the mean value when $\xd$ varies in $\T^d$, we obtain 
\begin{equation}
N \b_m = {\ob \over m} + l_m + o \left ( {1 \over m} \right ) .
\end{equation}

By $(3)$ and $(6)$, $l_m = r + o(1)$. As $l_m$ and $r$ are both vectors of $\Z^d$, $l_m = r$ if $m$ is large enough.
Equation $(6)$ becomes
\begin{equation}
N m \b_m = \ob + m r  + o (1).
\end{equation}
\noindent
Comparing $(2)$ and $(7)$, we have $k_m = m r + o(1)$ and this implies as above that $k_m = m r$ for large $m$. 
Equation $(2)$ now states that  
$$ \b_m = {\ob \over N m} + {r \over N}. $$

To finish the proof, simply define $i_m : \T^d \longrightarrow T^*\T^d$ as 
$$i_m = G \circ j_m = (\psi_m , f_m).$$
\noindent
Then $\psi_m(\xd) = \psi (\xd + u_m(\xd))$, so that the sequence converges in $C^\infty$ topology to 
$\psi$. The set $i_m(\T^d)$ is (as $j_m(\T^d)$) a Lagrangian manifold, and it is a graph because $\psi$ 
is a diffeomorphism.

\end{document}